\documentclass[twoside,11pt]{amsart}
\usepackage{amsmath,amssymb}

\newcommand{\CP}{\mathbb{CP}}
\newcommand{\bCP}{\overline{\mathbb{C P}}}

\newcommand{\TryPackage}[3]{\IfFileExists{#1.sty}{\usepackage{#1}#2}{#3}
}
\TryPackage{mathrsfs}{\renewcommand{\mathcal}{\mathscr}}{%
        \TryPackage{eucal}{}{}}

\newcommand{\ZZ}{{\mathbb Z}}
\newcommand{\RR}{{\mathbb R}}

\usepackage{latexsym}
\usepackage{amsmath}
\usepackage{amsthm}
\usepackage{amscd}
\usepackage{xypic} 
\usepackage{amsfonts}
\usepackage{amssymb}
\usepackage{graphicx}
\usepackage{psfrag}

\newtheorem{df}{Definition}
\newtheorem{thm}[df]{Theorem}
\newtheorem{scho}[df]{Scholium}

\newtheorem{lem}[df]{Lemma}
\newtheorem{prop}[df]{Proposition}

\begin{document}

\title{A symplectic manifold homeomorphic but not diffeomorphic to  $\CP^2 \# 3\overline{\CP}^2$}

\author{Scott Baldridge}
\author{Paul Kirk}
\date{February 2, 2007}

\thanks{The first   author  gratefully acknowledges support from the
NSF  grant DMS-0507857. The second  author  gratefully acknowledges
support from the NSF  grant DMS-0604310.}

\address{Department of Mathematics, Louisiana State University \newline
\hspace*{.375in} Baton Rouge, LA 70817}
\email{\rm{sbaldrid@math.lsu.edu}}

\address{Department of Mathematics, Indiana University \newline
\hspace*{.375in} Bloomington, IN 47405}
\email{\rm{pkirk@indiana.edu}}

\subjclass[2000]{Primary 57R17; Secondary 57M05, 54D05}
\keywords{Symplectic topology, Luttinger surgery, fundamental group, 4-manifold}

\begin{abstract}
In this paper we construct a minimal symplectic 4-manifold  and
prove it is homeomorphic but not diffeomorphic to $\CP^2\#
3\overline{\CP}^2$.
\end{abstract}

\maketitle


\section{Introduction}

The main result of this article is the construction of   a minimal symplectic 4-manifold  that is
homeomorphic but not diffeomorphic to $\CP^2\# 3\overline{\CP}^2$.

The construction of  manifolds homeomorphic but not diffeomorphic to
$\CP^2\# k\overline{\CP}^2$s for $k\leq 9$ began with Donaldson's
seminal example \cite{D} that the Dolgachev surface $E(1)_{2,3}$ is
not diffeomorphic to $ \CP^2\# 9\overline{\CP}^2$. In 1989, Dieter
Kotschick \cite{Kot} proved that the Barlow surface is homeomorphic
but not diffeomorphic to $\CP^2\# 8\overline{\CP}^2$.  In  2004
Jongil Park \cite{park1} constructed the first exotic smooth
structure on $\CP^2\# 7\overline{\CP}^2$.  Since then Park's results
have been expanded upon in \cite{OzS, SS, FS3, PSS}, producing
infinite families of smooth 4--manifolds homeomorphic but not
diffeomorphic to $\CP^2\# k\overline{\CP}^2$ for $k=5,6,7,8$.  The
$k=5$ examples  are not symplectic.

Akhmedov \cite{A} describes a construction of a symplectic
4-manifold homeomorphic to but not diffeomorphic to $\CP^2\#
5\overline{\CP}^2$.   Our approach is indebted to his  idea of using the symplectic sum
construction along genus 2 surfaces to kill fundamental groups in an efficient way. Earlier
approaches  start with a simply connected manifold and kill
generators of the second homology using the rational blowdown approach. \footnote{A similar result  to
our main theorem is announced in \cite{AP}.}

Using Luttinger surgery in addition to symplectic sums expands the
palette of available symplectic constructions, and combined with
Usher's theorem \cite{usher}, verifying that a construction yields a
minimal symplectic manifold is straightforward. This is the approach
taken in investigating small symplectic manifolds  in our
previous article \cite{BK2}, which among other things contains
examples of symplectic manifolds homeomorphic but not diffeomorphic
to  $\CP^2\# 5\overline{\CP}^2$.

Many of our constructions have their origin in \cite{FS4}, where
 symplectic sums of products of surfaces and surgery along nullhomologous tori
 are used to construct  symplectic and non-symplectic  manifolds  which are homeomorphic and in some cases not diffeomorphic.

Our experience, gleaned while working on \cite{BK, BK2, BK3},  taught us that
there are serious technical issues arising from working with
fundamental groups and cut and paste constructions, which can easily
lead to plausible but unverified or even incorrect calculations.  As
usual, base point issues are the culprit. Thus in writing the
present article we take great care in performing fundamental group
calculations. This is reflected in the length of the proof of
Theorem \ref{tough}, whose statement is perhaps not surprising in hindsight, but critical
for what follows.  At every stage of our constructions we must keep track not just of homotopy classes, but representative loops. We encourage the interested
reader to start with the proof of our main result, Theorem \ref{ex1}, and to  save the proof of Theorem 2 for last.

To summarize our construction, our example is the  symplectic sum of two manifolds along genus 2 surfaces. The first manifold $W$ is obtained  from Luttinger surgery on a pair of Lagrangian tori in
$T^4\#2\bCP^2$.  The second manifold $P$ is obtained by Luttinger surgery on    four  Lagrangian tori in $F_2\times T^2$, where $F_2$ is a surface of genus $2$.    Recall (\cite{Gompf}) that the symplectic sum is obtained by removing a neighborhood of a surface in each manifold, and gluing the resulting manifolds along their boundary.  Thus our approach is informed by the methods of knot theory: we essentially calculate the fundamental groups  of the complement of a link of two  tori and a genus 2 surface in $T^4\#2\bCP^2$ and the complement of a  link of four tori and a genus 2 surface in $F_2\times T^2$, as well as their meridians and longitudes with respect to paths from all the link components to the base point. It is this last point which makes the calculations challenging (and easy to screw up).

To make the exposition as concise as possible, we use the following strategy. To show a group is trivial, it suffices to show it is a quotient of the trivial group. More generally, one can view the Seifert-Van Kampen theorem as giving two pieces of information: first it provides generators and then identifies all relations. Since our goal is to  show that the example is simply connected, it suffices to find all generators and sufficiently many relations for the building blocks to reach the desired conclusion. Thus we eschew the problem of finding a complete presentation of the fundamental groups of $W$ and $P$, and content ourselves with establishing the relations we require for the proof.

We remark that the equation $\ell_2=bab^{-1}$ which appears in the statement of Theorem \ref{tough} (rather than the perhaps  expected $\ell_2=a$) hints at the fact  that calculations of fundamental groups of torus   surgeries on Lagrangian tori in the product  of surfaces  are likely to be subtle.   By stating Theorem \ref{tough} as we did (i.e.~in the product of punctured tori) it will be very useful in other contexts when small symplectic manifolds are to be constructed, since e.g.~ one can build products of closed surfaces starting with the product of punctured tori.

\section{fundamental group calculations}

Let $H$ be an oriented genus 1 surface with one boundary component. Let $x,y$   be oriented embedded circles  representing a symplectic basis of $H_1(H)$ so that $x$ and $y$ intersect transversally and positively in one point, which we denote by $h$.       Denote the corresponding based homotopy classes in  $\pi_1(H,h)$ also by $x$ and $y$

Now let $K$ be another  oriented genus 1 surface with one boundary component. Let $a,b$   be oriented embedded circles  representing a symplectic basis of $H_1(K)$ so that $a$ and $b$ intersect transversally and positively in one point, which we denote by $k$.

The image of the loops $x,y,a,b$ under the  inclusion $H\times\{k\}\cup \{h\}\times K\subset H\times K$  define homotopy classes which we as usual  denote by  $x,y,a,b\in \pi_1(H\times K, (h,k))$.      The base point $(h,k)$ for $H\times K$ is to be understood throughout this section.

\medskip

Let $X$ be a push off of $x$ in $H$ to the right with respect to  the orientations on $H$ and $x$.  Let  $Y$ be a parallel push off of $y$ to the left. Thus $x$ and $X$ are disjoint parallel curves on $H$.

Now let $A_1$ be a parallel push off of $a$ in $K$ to the right of $a$. Let  $A_2$ be a further parallel push off of $A_1$, to the right of $A_1$. Thus $a, A_1$ and $A_2$ are parallel curves in $K$.

Figure 1 illustrates all the curves on the surfaces $H$ and $K$.

\begin{figure}[h]
\bigskip
\begin{center}
\psfrag{y}{$y$} \psfrag{x}{$x$}
\psfrag{X}{$X$} \psfrag{h}{$h$}\psfrag{K}{$K$}\psfrag{H}{$H$}
\psfrag{Y}{$Y$} \psfrag{A1}{$A_1$} \psfrag{A2}{$A_2$} \psfrag{a}{$a$} \psfrag{b}{$b$} \psfrag{k}{$k$}
  \includegraphics[scale=1]{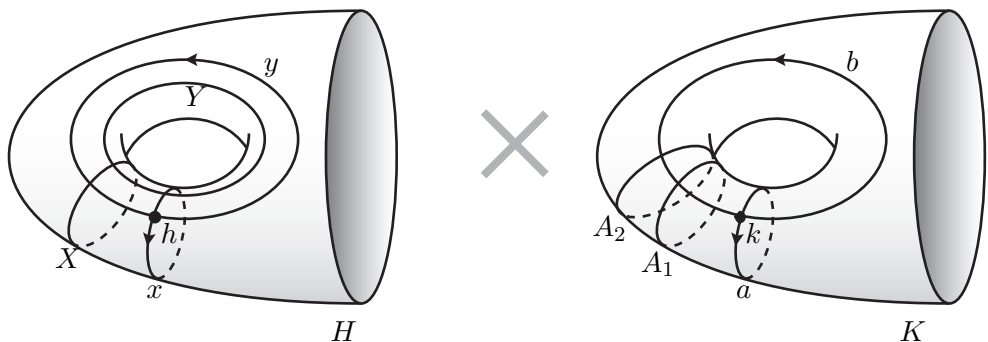}
  \caption{The surface $H\times K$.}
  \label{fig2}
\end{center}
\end{figure}

\bigskip

We define two disjoint tori $T_1,T_2$ in $H\times K$ as follows.

$$T_1= X\times A_1   \text{ and } T_2= Y \times A_2.$$

Fix a product symplectic form on $H\times K$. (Typically we think of $H$ and $K$ as codimension 0 submanifolds of   closed tori $\hat{H}$ and $\hat{K}$ and restrict the standard  product symplectic form on $\hat{H}\times \hat{K}$ to $H\times K$.)

\medskip

 The proof of the following proposition is simple.

 \begin{prop}  The tori $T_1 $ and $T_2$  are Lagrangian and the  surfaces $H\times \{k\}$ and $\{h\}\times K$ are symplectic.  Moreover, $T_1$ and $T_2$ are  disjoint and disjoint from
 $H\times \{k\}$ and $\{h\}\times K$.   \qed

\end{prop}

Notice that every torus of the form $C\times D\subset H\times K$,
(where $C\subset H$ and $D\subset K$ are embedded curves) is
Lagrangian.  Recall that a Lagrangian torus $T$ in a symplectic
4-manifold $M$ has a canonical framing called the {\em  Lagrangian
framing}. In fact, the Darboux-Weinstein theorem  \cite{DS}
implies
that
 a tubular neighborhood of $T$ can be identified with $T\times D^2$ in such a way that the parallel tori in $M$ corresponding to $T\times \{d\}$ in this framing are also Lagrangian for every $d\in D^2$.   In particular, given any such neighborhood and any $d\in \partial D^2$, we will call the torus $T\times \{d\}$ in the boundary of a tubular neighborhood of $T$ a {\em Lagrangian push off} of $T$, and if $\gamma\subset T$ is a curve we call the curve corresponding $\gamma\times \{d\}$ the {\em Lagrangian push off of $\gamma$}.

\medskip

The following theorem is the critical step in our constructions. Before we state it, we begin with an  observation and a warning. First the observation: the torus $T_2$ intersects the torus $x\times b$ transversally in one point. Together with the remarks about the Lagrangian framing discussed above, one concludes without much trouble that
in $\pi_1(H\times K-(T_1\cup T_2))$,  the meridian of $T_2$ takes the form $[\tilde{x},\tilde{b}]=
\tilde{x}\tilde{b}\tilde{x}^{-1}\tilde{b}^{-1}$, and the Lagrangian push off of the curves $Y$ and $A_2$ take the form $\tilde{y}$ and $\tilde{a}$ respectively, where for $z\in\pi_1(H\times K-(\cup_iT_i))$ we let $\tilde{z}$ denote  some conjugate of $z$.

Put another way, consider the three circles  that lie on the boundary of a tubular neighborhood of $T_2$, namely the boundary of a meridian disk $\{t\}\times D^2$, and the Lagrangian push offs of the  curves $Y$ and $A_2$ with respect to a normal Lagrangian vector field. These curves are freely homotopic to (respectively) the triple $[x,b], y$ and $a$ in $H\times K-(T_1\cup T_2)$.

But {\em they need not be equal to this triple}  when the boundary of the tubular neighborhood is joined by a path to the base point $(h,k)$ in $H\times K-(T_1\cup T_2)$.   There is some freedom in the choice of path to simultaneously conjugate all three. But
to expect that there exists a path to the base point so that $( [\tilde{x},\tilde{b}],\tilde{y},\tilde{a})=([x,b],y,a)$ in $\pi_1(H\times K-(T_1\cup T_2))$   is in general too much to hope for, and has led to some confusion and mistakes which we need to avoid.

 The configuration  is nevertheless sufficiently   explicit in our situation to prove the following theorem.

\begin{thm} \label{tough} There exist paths in    $H\times K-(T_1\cup T_2)$ from the base point to the boundary of the tubular neighborhoods $T_1\times \partial D^2$ and $T_2\times \partial D^2$ with the following property. Denote  by  $\mu_i, m_i, \ell_i\in \pi_1(H\times K-(T_1\cup T_2))$ the  loops obtained by following the path to the boundary of the tubular neighborhood of $T_i$, then following (respectively) the meridian of $T_i$ and the two  Lagrangian push offs of the generators  on $T_i$ are given by the following formulae:
$$\mu_1=[ b^{-1}, y^{-1}], m_1=x, \ell_1= a, $$
and
$$\mu_2=[x^{-1},b], m_2=y, \ell_2=bab^{-1}.$$
where $x,y,a,b\in \pi_1(H\times K -(T_1\cup T_2))$ are the loops which lie on the surfaces $H\times \{k\}$ and $\{h\}\times K$ described above.

Moreover,
$\pi_1(H\times K -(T_1\cup T_2))$ is generated by $x,y,a,b$ and the relations
$$[x,a]=1, [y,a]=1, [y,bab^{-1}]=1$$
as well as
$$[[x,y],b]=1,[x,[a,b]]=1,[y,[a,b]]=1$$
hold in $\pi_1(H\times K -(T_1\cup T_2)).$
\end{thm}

\noindent{\bf Remark.} Note that we are not assuming any particular orientation convention on the meridians, or even that the two meridians are oriented by the same convention.  Looking ahead, when we perform Luttinger surgeries below we are free to do either $ 1$  or $-1$ surgeries, and we will pick the sign that introduces the relation we require.

\begin{proof}   Figures 2 and 3 will guide the reader through the argument. View a torus as a square with opposite sides identified, thus $T^4$ can be thought of as a quotient of the product of two squares. Equivalently, we think of it as a quotient of the cube with coordinates $x,y,b$ and an interval corresponding to the $a$ coordinate.   Since $H\times K\subset T^4$, we visualize $H\times K$ as a subset of the 4-cube.

\bigskip

\begin{figure}[h] 
\begin{center}
\bigskip
\psfrag{y}{$y$} \psfrag{x}{$x$} \psfrag{X}{$X$}
\psfrag{h}{$h$}\psfrag{K}{$K$}\psfrag{H}{$H$}\psfrag{D}{$X \times
A_1$}\psfrag{YA2}{$Y\times A_2$} \psfrag{Y}{$Y$} \psfrag{A1}{$A_1$}
\psfrag{A2}{$A_2$} \psfrag{a}{$a$} \psfrag{b}{$b$}
\psfrag{k}{$k$}\psfrag{hk}{$(h,k)$}\psfrag{bdH}{$\partial H$}
\includegraphics[scale= .85]{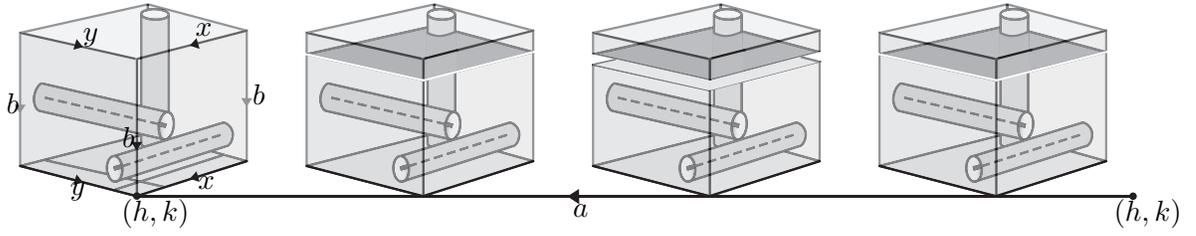}
 \caption{$H\times K$}
  \end{center}
  \bigskip
 \end{figure}

\bigskip

 \begin{figure}[h] 
 \begin{center}
 \bigskip
\psfrag{y}{$y$} \psfrag{x}{$x$} \psfrag{X}{$X$}
\psfrag{h}{$h$}\psfrag{K}{$K$}\psfrag{H}{$H$}\psfrag{D}{$X \times
A_1$}\psfrag{YA2}{$Y\times A_2$} \psfrag{Y}{$Y$} \psfrag{A1}{$A_1$}
\psfrag{A2}{$A_2$} \psfrag{a}{$a$} \psfrag{b}{$b$}
\psfrag{k}{$k$}\psfrag{hk}{$(h,k)$}\psfrag{bdH}{$\partial H$}
\includegraphics[scale= 1.1]{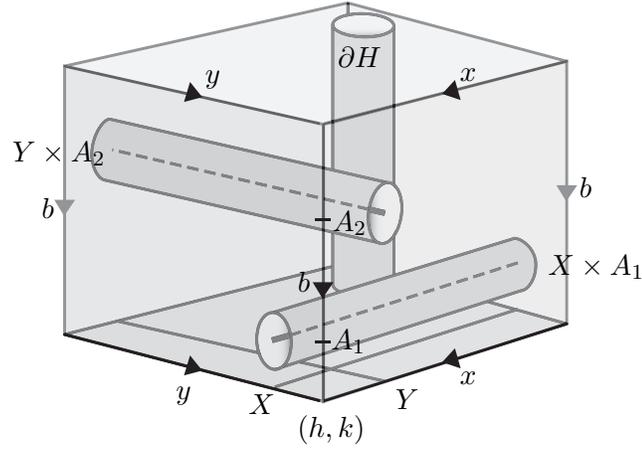}
 \caption{The slice $a=1$}
 \bigskip
 \end{center}
 \end{figure}

 We start with the easy torus $T_1$ first.  Let $p\in H$ denote the intersection point of $X$ and $y$. Let $\alpha$ be the following path from the base point to the boundary of the tubular neighborhood of $T_1$.  Starting at $(h,k)$, let $\alpha_1$ denote the path traced out by traveling backwards along $y$ in $H\times \{k\}$ until you hit $X$ at the point $(p,k)$. Then  let $\alpha_2$ denote the path obtained by traveling  in $
\{p\}\times K$ backwards along $b$ until just before you hit $A_1$. This defines the path $\alpha=\alpha_1*\alpha_2$ in $H\times K-(T_1\cup T_2)$ from the base point to the point $(p,q)$, where $\{p\}=X\cap y$  and $q$ is a point on $b$ just to the right of $A_1$.

The square $[0,1]\times [0,1]$ maps to $H\times K-  T_2 $ by
$$  y^{-1}\times b^{-1}:I\times I\to H\times K-  T_2 .$$
The interior of this square intersects $T_1$ transversally once.
Moreover, the path $\alpha$  lies on this square starting at the image of  $(0,0)$. It follows that the meridian of $T_1$ is (based) homotopic to the boundary of this square, starting at $(0,0)$, i.e. $\mu_1=[b^{-1},y^{-1}]$.

Next consider the loop $m_1$ which follows $\alpha$ to $(p,q)$, then follows the loop $X\times \{q\}$  around back to $(p,q)$, and finally returns to the base point along $\alpha^{-1}$. This is the Lagrangian push off of $X$ since the second coordinate $q$ is held fixed as one moves along $X$.

We show that the loop $m_1$ is based homotopic to $x$ in $H\times K-(T_1\cup T_2)$.  First, there is an annulus   in $H$ with boundary $x$ and $X$ which contains the
arc   from $h$ to $p$ following $y$ backwards. This determines an annulus $F_1$ in $H\times\{k\}\subset H\times K$ which misses $T_1\cup T_2$ with $(h,k)$ on one boundary circle, $(p,k)$ on the other, and the arc $\alpha_1$ spanning these two points. There is another annulus $F_2$ of the form $X\times \alpha_2$ which contains the arc $\alpha_2$ and misses $T_1\cup T_2$. Gluing $F_1$ to $F_2$ along their common boundary $X\times \{k\}$ yields a homotopy from $x$ to $m_1$  which is base point preserving since it contains the path $\alpha$ spanning the two boundary components.

Next, consider the loop $\ell_1$ which first follows $\alpha$ to $(p,q)$, then follows the loop $\{p\}\times A_1^+$ where $A_1^+$ is the parallel copy of $A_1$ in $K$ that passes through $q$, and finally returns to the base point along $\alpha^{-1}$.  As explained above, $\{p\}\times A_1^+$ is the Lagrangian push off of $A_1\subset T_1$ since it is the push off  of $A_1$ in $K$.

We show that the loop $\ell_1$ is based homotopic to $a$.  We argue similarly as above. This time there is an annulus $F_3$ which lies in $H\times a$ with boundary the curves $\{h\}\times a$ and  $\{p\}\times a$ which contains the path $\alpha_1$ spanning its boundary components.  There is an annulus $F_4$ in $\{p\}\times K$ with boundary the curves $\{p\}\times a$ and $\{p\}\times A_1^+$ which contains the path $\alpha_2$. This proves that $a$ and $\ell_1$ are based homotopic.

We now turn to the other torus $T_2$. The attentive reader will realize that the difficulty here is that the analogue of the path $\alpha_2$ we would want to use intersects $T_1$. The solution presents itself from this consideration: we will need to travel {\em forwards} along $b$ until we approach $A_2$

Proceeding in earnest now, let $r\in H$ denote  a point on $x$ close to and to the right of $Y$ (and left of $y$.) Let $s\in K$ denote the intersection point of $A_2$ with $b$.  Let $\beta_1$ be the path in $\{h\}\times K$ which starts at $(h,k)$ and moves  forwards  along $\{h\}\times b$ to the point $(h,s)$.  Let $\beta_2$ be the path in $H\times \{s\}$ starting at $(h,s)$ and moving along $x$ backwards until  the point $(r,s)$ in the boundary of the tubular neighborhood of $T_2$ is reached.
 The path $\beta=\beta_1*\beta_2$ is our path from the base point to the boundary of the tubular neighborhood of $T_2$.

To compute $\mu_2$, we notice that there is a map of a square:
$$x^{-1}\times b: I\times I\to H\times K- T_1$$
which intersects $T_2$ transversely once and contains the path $\beta$, starting at $(0,0)$. Thus $\mu_2$ can be read off the boundary of the square, and hence $\mu_2=[x^{-1},b]$.

Next, consider the loop $m_2$ which follows $\beta$ to $(r,s)$, then follows $Y^+\times \{s\}$ and returns to the base point along $\beta^{-1}$, where $Y^+$ is the push off of $Y$ in $H$ which passes through $r$.  This is the Lagrangian push off of $Y$, since $Y^+\times\{s\}$ is a Lagrangian curve.  There is an annulus $F_5$  with boundary $y\times \{h\}$ and $y\times \{s\}$ which contains the path $\beta_1$. There is   an annulus $F_6\times\{s\}$ with boundary $y\times\{s\}$ and $Y^+\times\{s\}$ which contains the path $\beta_2$. These glue to give a base point preserving homotopy of   $m_2$ to $y$.

We saved the most difficult calculation for last, and it is here that Figure 2 becomes most helpful. Consider the loop $\ell_2$ which follows $\beta$ to $(r,s)$, then follows $\{r\}\times A_2$ and then returns along $\beta^{-1}$. There is a  surface $F_7$ in
$\{h\}\times K$ (a punctured annulus) with three boundary components: $\{h\}\times a$, $\{h\}\times A_2$, and $\{h\}\times \partial K$ which contains the path $\beta_1$.  There is an annulus $F_8$ of the form $\beta_2\times A_2$ with boundary $\{h\}\times A_2$ and $\{r\}\times A_2=\ell_2$.

Cut a slit in $F_7$ along an arc of the form $\{h\}\times \gamma$, where $\gamma$ is a path in $K$ from $k$ to the boundary.   Then the commutator $bab^{-1}a^{-1}$ is homotopic to the composite of $\gamma$, the loop that follows the boundary, and then $\gamma^{-1}$.  Cutting $F_7$ along $\beta_1$ and $\gamma$ and reading the word on the boundary one finds
$\beta_1*A_2^{-1}*\beta_1^{-1}* bab^{-1}a^{-1}*a$ and gluing on $F_8$ one concludes that
$$\ell_2=bab^{-1}a^{-1}a= bab^{-1}.$$

(For the benefit of the reader, we sketch an alternative way to see this, referring to Figure 3. Let $\beta_3$ be the path following $b$ forwards starting at $\beta_1(1)$, so that $\beta_1*\beta_3=b$. The square   of the form $\beta_2^{-1}\times a$ glues to the square $\beta_3\times a$ to give a homotopy from $\ell_2$ to $bab^{-1}$.)

\medskip

We now turn to the assertions about $\pi_1(H\times K-(T_1\cup T_2))$.  The surface $K$ decomposes into two  surfaces: an annulus $K_1$ with boundary $A_1$ and $A_2$ and its complement, a 3-punctured sphere with boundary the disjoint union $\partial K \cup A_1\cup A_2$.

 We take the preimages of the $K_i$ via the projection to $K$. Precisely, let $\Phi:H\times K\to K$ denote the projection and define
$$W_1=  \Phi^{-1}(K_1) \cap (H\times K- nbd(T_1\cup T_2))$$
and
$$W_2=  \Phi^{-1}(K_2) \cap (H\times K- nbd(T_1\cup T_2))$$
Notice that $W_1$ is homeomorphic to $H\times K_1$ and $W_2$ is homeomorphic to $H\times K_2$.

Thus $W_1\cup W_2=H\times K- nbd(T_1\cup T_2)$.  The intersection $W_1\cap W_2$ has two  components: one of them is
$$ \Phi^{-1}(A_1) \cap (H\times K- nbd(T_1\cup T_2))= H\times A_1-nbd(T_1)=(H-nbd(X))\times A_1.$$
The other one is
$$ \Phi^{-1}(A_2) \cap (H\times K- nbd(T_1\cup T_2))= H\times A_2-nbd(T_2)=(H-nbd(Y))\times A_2.$$

\medskip

To apply the Seifert-Van Kampen theorem requires the intersection to be connected, so we take the usual approach (e.g. taken when computing fundamental groups of bundles over $S^1$) and change $W_1$ and $W_2$ slightly to make their intersection connected as follows.

Let $\tau$ denote the arc in $\{h\}\times K$ which starts at the base point $(h,k)$ and travels along $b$ {\em backwards}, passes through $A_1$, and ends at the intersection point of $b$ with $A_2$.

We let $W_1'=W_1\cup \tau$, this is just $W_1$ with a small hair attached connecting it to the base point.    Define three loops in $W_1'$ based at $(h,k)$ as follows. Let $k'$ denote a point on $b$ between $A_1$ and $A_2$.   Follow  the arc $\tau$ from $(h,k)$ to $(h,k')$, then take the loop $x\times \{k'\}$, then return to $(h,k)$ along $\tau^{-1}$. Call this loop $x'$. Similarly define the loop $y'$. Finally, define the loop $a'$ to be the loop obtained by following $\tau$ from $(h,k)$ to $(h,k')$, then following a loop parallel to and between $\{h\}\times A_1$ and $\{h\}\times A_2$ in $\{h\}\times K$, and finally returning to $(h,k)$ along $\tau^{-1}$.

 Since $W_1$ is homeomorphic to $H\times K_1$ (always taking the base point $(h,k)$),
$$\pi_1(W_1')=\langle x', y'\rangle \oplus \ZZ a'.$$

We then let $W_2'=W_2\cup \tau$.  This is $W_2$ with an arc attached spanning the boundary components corresponding to $A_1$ and $A_2$.   Notice that  the loops $a,b,x,$ and $y$ all lie  in $W_2'$ (recall that these are the explicit loops on $H\times \{k\}\cup \{h\}\times K$ which we claim generate $\pi_1(H\times K-(T_1\cup T_2))$).  Denote by $c$ the loop in $K_2$ based at $k$ which
travels to the boundary $\partial K$, goes around once, and returns to $k$ (thus, in $\pi_1(K)$, $c$ represents the commutator $[a,b]$). We consider the loop $c':=\{h\}\times c$ in $\{h\}\times K$ based at $(h,k)$, this is also a loop in $W_2'$.

Then because  $W_2'$ is obtained from $W_2\cong H\times K_2$ by adding the arc $\tau$, it is clear that the five loops $a,b,x,y,c'$ generate $\pi_1(W_2')$.  We will not need this, but note that $a$ and $c'$ commute with $x,y$ and that $b$ generates a free factor.

We now apply the Seifert-Van Kampen theorem to conclude that
 $\pi_1(H\times K- nbd(T_1\cup T_2))$ is generated by the loops
 $$a',x',y', a,b,x,y,c'.$$
 Thus to establish our claim that  $x,y,a,b$ generate $\pi_1(H\times K- nbd(T_1\cup T_2))$,
we must show that the based homotopy classes $a',x',y',$ and $c'$  in
$\pi_1(H\times K- nbd(T_1\cup T_2))$ can be expressed in terms of $a,b,x,$ and $y$.

Since $a'$ lies on $\{h\}\times K$, which misses $T_1\cup T_2$, it is obvious that $a'$ and $a$ represent the same class. Equally easy is the observation that $c'=[a,b]$ in $\pi_1(H\times K- nbd(T_1\cup T_2))$.

This leaves the classes $x'$ and $y'$.  Consider first $x'$. We claim it is based homotopic to $x$.
We can give an explicit formula for such a homotopy. Let $\beta$ denote the path from $k$ to $k'$ in $K$ that follows $b$ backwards. For $s\in [0,1]$, let  $\beta_s$ denote the path $t\mapsto\beta((1-s)t)$ (so $\beta_0=\beta$ and $\beta_1$ is the constant path at $k$).

   Then the homotopy
$$s\mapsto (\{h\}\times \beta_s)*(x\times \{\beta(1-s)\})* (\{h\}\times \beta_s)^{-1}$$
is a based homotopy  from $x'$ to $x$ {\em that misses $T_1\cup T_2$}.  This is because ,when passing through $ \Phi^{-1}(A_1) \cap (H\times K- nbd(T_1\cup T_2)) =(H-nbd(X))\times A_1$ (i.e. when $\beta_s(1)=\beta(s)$ lies on $A_1$), the curve $x$ is parallel to $X$ and hence misses it.

We can similarly show that $y'$ is based homotopic in $H\times K- nbd(T_1\cup T_2)$ to a loop represented by a word in $a,b,x,y$. This time we need to push across $A_2$ instead of $A_1$. Since we have already noticed that any based loop in $W_2'$ can be expressed in terms of $a,b,x,y$, it is easiest just to slide $y'$ along $\tau$ past $A_2$. This time the fact that
$Y$ is parallel to $y$ and $ \Phi^{-1}(A_2) \cap (H\times K- nbd(T_1\cup T_2))=(H-Y)\times A_2$
allows us to conclude that $y'$ can be expressed in terms of $a,b,x,y,c'$, and hence $a,b,x,$ and $y$.

Thus $\pi_1(H\times K- nbd(T_1\cup T_2))$ is generated by the loops $a,b,x,y$, as claimed.

 To finish the proof, we establish the stated commutator relations. The torus $x\times a$ contains the base point $(h,k)$ and the curves $x$ and $a$, and misses $T_1$ and $T_2$. Hence $[x,a]=1$. Similarly the torus $y\times a$ shows that $[y,a]=1$. If $e$ denotes the loop in $H$ that goes from the base point to the boundary of $H$,   travels around the boundary, then returns to $h$ (avoiding the curves $X$ and $Y$) then the mapped in torus $e\times b$ misses $T_1\cup T_2$ and hence $[[x,y],b]=1$ in $\pi_1(H\times K- nbd(T_1\cup T_2))$. Similarly $[x,[a,b]]=1$ and $[y,[a,b]]=1$.

   The other  commutator relation is a consequence of the fact that $\mu_i,m_i,$ and $\ell_i$ commute, since they live on the boundary of a tubular neighborhood of $T_i$, a 3-torus.

This completes the proof of Theorem \ref{tough}.
\end{proof}

\medskip

There are a few other  relations in  $\pi_1(H\times K- nbd(T_1\cup T_2))$ which we did not mention in the statement of Theorem \ref{tough}, e.g. $[[b^{-1},y^{-1}],x]$, $[[b^{-1},y^{-1}],a]$, $[[x^{-1},b],y]$,  and $[[x^{-1},b], bab^{-1}]$.  These follow from the fact that they correspond to loops on the boundary of the tubular neighborhoods of $T_1$ and $T_2$. We will not need these relations in our argument.

\bigskip

In Theorem \ref{tough} we worked with the product of two punctured tori not for generality's sake, but because we will need to use the same construction in three different contexts later:
\begin{enumerate}
\item $H$ is the complement of a disk in a (closed) torus $\hat{H}$. Thus we will be interested in the
two Lagrangian tori $T_1, T_2$ in $\hat{H}\times K$,  the fundamental group
$\pi_1(\hat{H}\times K-(T_1\cup T_2))$, and the corresponding $\mu_i\,m_i,\ell_i$.
\item $K$ is the complement of a disk in a (closed) torus $\hat{K}$. Thus we will be interested in the
two Lagrangian tori $T_1, T_2$ in $H\times \hat{K}$, the fundamental group
$\pi_1(H\times \hat{K}-(T_1\cup T_2))$, and the corresponding $\mu_i\,m_i,\ell_i$.
\item $K$  and $H$ are both  complement of  disks in a (closed) tori. Thus we will be interested in the
two Lagrangian tori $T_1, T_2$ in the four torus $T^4$, the fundamental group
$\pi_1(T^4 -(T_1\cup T_2))$, and the corresponding $\mu_i\,m_i,\ell_i$.
\end{enumerate}

(Cases (1) and   (2) are inequivalent due to the asymmetry of the pair $X,Y$ and the pair $A_1,A_2$). The effect on fundamental groups in these three cases is clearly to impose the appropriate commutator relation.

\begin{scho} \label{scho}  In the three cases enumerated above, the statement of Theorem \ref{tough}
remains true if we replace $H\times K$ by $\hat{H}\times K,  H\times \hat{K}$, and $T^4$ respectively.   Moreover, in the three cases, there is a further relation in the  fundamental groups:
\begin{enumerate}
\item The relation $[x,y]=1$ holds in
$\pi_1(\hat{H}\times K-(T_1\cup T_2))$.
\item The relation $[a,b]=1$ holds in
$\pi_1(H\times \hat{K}-(T_1\cup T_2))$.
\item The relations $[x,y]=1$  and $[a,b]=1$ hold  in
$\pi_1(T^4 -(T_1\cup T_2))$.
\end{enumerate}
\qed
\end{scho}

\bigskip

Recall that given a Lagrangian torus $T$ in a symplectic 4-manifold $M$, with meridian $\mu$, and Lagrangian push offs $m$ and $\ell$ in $\pi_1(M-T)$, {\em Luttinger surgery} is the process which removes a neighborhood $T\times D^2$ from $M$  and glues it back in by a diffeomorphism which takes
a disk $\{t\}\times D^2$ to a curve of the form $\mu m^{kp}\ell^{kq}$ where $p,q$ are relatively prime integers and $k$ is an integer.  To specify the choices, we say the resulting manifold is obtained
{\em by $1/k$ Luttinger surgery along the curve $pm+q\ell$.}  Luttinger \cite{Lut} (see also \cite{ADK}) proved that for any integer $k$ and any choice of $p,q$, the result of Luttinger surgery admits a symplectic structure in which the core $T\times \{0\}$ is also Lagrangian, and so that the symplectic structure is unchanged in the complement of the tubular neighborhood of $T$.

We include the following well-known lemma for completeness.

\begin{lem} \label{lut} The fundamental group of the manifold obtained by $1/k$  Luttinger surgery on $M$ along $pm+q\ell$ is the quotient
$$\pi_1(M-T)/N(\mu m^{kp}\ell^{kq})$$
where $N(\mu m^{kp}\ell^{kq})$ denotes the normal subgroup of $\pi_1(M-T)$ generated by
$\mu m^{kp}\ell^{kq}$
\end{lem}
\begin{proof}  The 2 torus has a handle structure with one 0-handle, two 1-handles, and one 2-handle. Thus the product $T^2\times D^2$ has a handle structure with one 0-handle, two 1-handles, and one 2-handle.  Looking from the outside in, one sees that attaching $T^2\times D^2$ can be accomplished by attaching one 2-handle,  two 3-handles, and one 4-handle.  Attaching the  2-handle has the stated effect on fundamental groups, and attaching 3 and 4 handles does not further affect the fundamental group.
\end{proof}

Call the relations in Theorem \ref{tough} and Scholium \ref{scho} {\em universal relations} since they hold for any Luttinger surgery, and indeed, in the complement of $T_1\cup T_2$. The relations of  Lemma \ref{lut} coming from Luttinger surgery  will be called {\em Luttinger relations}.

\bigskip

We end this section  with one lemma which will be used to establish minimality of the manifolds we construct.

\begin{lem}\label{mini}
Let $M$ be obtained from the 4-torus $T^4=\hat{H}\times\hat{K}$ by $1/k_1$ Luttinger surgery on $T_1$ along $x$ and $1/k_2$ surgery on $T_2$ along $a$. Then $\pi_2(M)=0$, and hence $M$ is minimal.
\end{lem}
\begin{proof}

First,  $1/k_1$ surgery on $T_1$ along $x$ transforms $T^4$ into $N\times S^1$, where $N$ is the 3-manifold that fibers over $S^1$ with monodromy  the $k_1$th power of the Dehn twist  on $\hat{H}$ along $x$.    This follows from the well-known fact for fibered 3-manifolds that changing the monodromy by a Dehn twist corresponds to a Dehn surgery along a curve in a fiber. One can find a careful  explanation  in \cite[pg. 189]{ADK}.

View $N\times S^1$ as a trivial circle bundle over $N$.   Removing a neighborhood of $T_2$ and regluing has the effect of changing this trivial $S^1$ bundle to a non-trivial bundle. Explicitly one removes a neighborhood of $y$ in $N$ and its preimage in $N\times S^1$, then reglues in such a way that $k_2[y]$ becomes  the divisor of the resulting $S^1$ bundle. Details can be found in \cite{B1}. In any case one can check directly from the construction that $M$ has a free circle action which coincides with the action on $N\times S^1$ away from $T_2$.

Thus $M$ is an $S^1$ bundle over a fibered 3-manifold $N$ with fiber a torus. It follows from the long exact sequence of homotopy groups that $\pi_2(M)=0$,  and hence $M$ contains no essential 2-spheres. In particular, $M$ is minimal.
\end{proof}

\section{The building blocks}

 \subsection{The manifold $W$}

 Consider the 4-torus $T^4=S^1\times S^1\times S^1\times S^1=T^2\times T^2$. Denote the coordinate  circles respectively by
 $s_1,t_1,s_2,t_2$. So for example $s_2=\{1\}\times \{1\}\times S^1\times \{1\}$.  These determine loops in $T^4$.  Let $\Phi:T^4\cong\hat{H}\times \hat{K}$ be a base point preserving diffeomorphism (in fact linear map) that takes the
 circles $s_1,t_1,s_2,t_2$ to $ x,y,a,b$ respectively.  Pulling back the tori $T_1,T_2$ via the symplectomorphism $\Phi$ gives a pair of Lagrangian tori in $T^4$ which, by abuse of notation, we also denote by $T_1$ and $T_2$.

It is helpful to  call $T_1$ the {\em  $s_1\times s_2$ torus} and $T_2$ the {\em  $t_1\times s_2$ torus}  to remember what (conjugacy) classes in the fundamental group they carry. The nomenclature can be confusing, since $T_2$ is pushed farther away that $T_1$ from the loop $a$, due to the fact that $A_1$ and $A_2$ are different curves in $K$. In particular, the Lagrangian push offs are only specified up to conjugacy by this notation: for example, Theorem \ref{tough} states that the Lagrangian push off of the curve on $T_2$ represented by $s_2$ curve is $\ell_2=t_2s_2t_2^{-1}$.

Theorem \ref{tough}, Scholium \ref{scho}, and Lemma \ref{lut} allows us to conclude that the fundamental group of the manifold $V$ obtained by $-1$ Luttinger surgery on the  $s_1\times s_2$ torus   along $s_1$ and $-1$ surgery on the $t_1\times s_2$ torus along  $s_2$ is generated by
  $s_1,t_1,s_2,t_2$
 and the Luttinger relations
$$[t_2^{-1},t_1^{-1}]=s_1 , [s_1^{-1},t_2]=t_2s_2t_2^{-1}$$
 as well as the universal relations
$$[s_1,t_1]=1, [s_2,t_2]=1, [s_1,s_2]=1, [t_1,s_2]=1, [t_1,t_2s_2t_2^{-1}]=1
$$
  hold. Note that   by conjugating by $t_2^{-1}$ we may  simplify the second Luttinger relation to
  $$[t_2^{-1},s_1^{-1}]=s_2.$$
 The last universal relation  reduces to the (redundant) relation $[t_1,s_2]=1.$

 Thus $\pi_1(V)$ is a quotient  of  the group with presentation  $$\langle s_1,t_1,s_2,t_2\ | \ [s_1,t_1], \ [s_2,t_2],\  [s_1,s_2], \ [t_1,s_2],
 \ [t_2^{-1},t_1^{-1}]s_1^{-1},\  [t_2^{-1},s_1^{-1}]s_2^{-1} \rangle.$$

\noindent{\bf Remark.} It is critical in these calculations   that the loops $s_1,t_1$ are to be understood as explicit loops in the symplectic surface $\hat{H}=T^2\times\{(1,1)\}\subset  T^4-(T_1\cup T_2)$ and the loops $s_2,t_2$ are to be understood as   loops in the symplectic surface
$\hat{K}=\{(1,1)\}\times T^2\subset  T^4-(T_1\cup T_2)$, all based at $(1,1,1,1)$.

\medskip

Lemma \ref{mini} shows that  $V$ can be described   as an $S^1$-bundle over
a 3-manifold that fibers over a circle with genus one fibers, and  so $V$ is a minimal symplectic 4-manifold.

\bigskip

The  symplectic tori  $\hat{H} =T^2\times\{(1,1)\}$ and $ \hat{K}=\{(1,1)\}\times T^2$ in  $T^4$ miss
neighborhoods of $T_1$ and $T_2$, and hence determine symplectic tori in $V$ that we continue to call
$\hat{H}$ and $ \hat{K}$.  Notice that  $\hat{H}$ and $ \hat{K}$ intersect  once transversally and positively at the base point $p=(h,k)=(1,1,1,1)$.

Symplectically resolve this intersection point as explained in \cite{Gompf}. This is a local modification in a small neighborhood of $p$ which replaces $\hat{H}\cup \hat{K}$ by a smooth symplectic surface $G$.

The topological description of this process is as follows. In a small 4-ball around $p$,  a pair of intersecting 2-disks in $\hat{H}\cup \hat{K}$ are removed and replaced by an annulus so that the resulting closed genus 2 surface $G$ is oriented consistently with the orientations of $\hat{H}$ and $\hat{K}$.   Thus one can choose a base point $p'$ inside this annulus, based loops  $s_1',t_1',s_2',t_2'$  on $G$ satisfying $[s_1,t_1][s_2,t_2]=1$ in $\pi_1(G,p')$,
and a small arc in the 4-ball from $p'$ to $p$ so that the inclusion $\pi_1(G,p')\to \pi_1(V,p')$ followed by the identification $\pi_1(V,p')\cong\pi_1(V,p)$ given by the small arc takes $s_i',t_i'$ to $s_i,t_i$.   Therefore we can safely rename $p'=p, s_i'=s_i$, $t_i'=t_i$ and the fundamental group calculations are unchanged.

\medskip

Now blow up $V$    twice at two distinct points on $G$, obtaining a symplectic manifold
$$W=V\#2\bCP^2.$$
 The proper transform of $G$ is a symplectic surface in $W$ (\cite{Gompf})  which we continue to call $G$.   It has  the same fundamental group properties as it did in $V$, but, in addition, $G\subset W$ has a trivial normal bundle and intersects each exceptional sphere transversally once.

 Fix a push off $G\to W-nbd(G)$ and give $W-G$ the base point which is the image of $p$ via this push off. Use a path in a meridian disk to identify based loops in  $W-G$  and based loops in $W$.   Since the surface $G$ intersects a sphere (either of
the two exceptional spheres) transversally in one point, the meridian of $G$ in $W-G$ is
nullhomotopic. Moreover, the inclusion  $W-G \subset W$ induces
an isomorphism  on fundamental groups, since every loop in $W$ can
be pushed off $G$ and every homotopy that  intersects $G$ can be
replaced by a homotopy that misses $G$ (using the exceptional sphere and the fact that $G$ is connected).  Therefore we conclude the following lemma. As before $N(S)$ denotes the normal subgroup generated by a set $S$.

\begin{lem}\label{lem1} \hfill
\begin{enumerate}
\item  The closed symplectic 4-manifold $W$ contains a closed symplectic genus 2 surface $G$ with trivial normal bundle.   There are  based loops $s_1,t_1,s_2,t_2$  on $G$ representing a standard symplectic generating set for $\pi_1(G,p)$ (thus satisfying $[s_1,t_1][s_2,t_2]=1$) such that these  loops generate $\pi_1(W,p)$ and, in
  $\pi_1(W,p)$, the relations
$$ 1=
[s_1,t_1]=[s_2,t_2]= [s_1,s_2]=[t_1,s_2]= [t_2^{-1},t_1^{-1}]s_1^{-1}=  [t_2^{-1},s_1^{-1}]s_2^{-1} $$
hold.  The inclusion $W-G\subset W$ induces an isomorphism on fundamental groups.

\item
Let $R$ be any 4-manifold containing a genus 2 surface $F$ with
trivialized normal bundle.  Let $\phi:G\to F$ be a diffeomorphism,
and set  $g_i=\phi_*(s_i),h_i=\phi_*(t_i)$ in $\pi_1(R)$.  Given a
map  $\tau:G\to S^1$, let  $\tilde{\phi}:G\times S^1\to F\times S^1$
the diffeomorphism given by $\tilde{\phi}(a,s)=(\phi(a),
\tau(a)\cdot s)$.
Form the sum:
$$S = (R-nbd(F))\cup_{\tilde{\phi}}(W-nbd(G)).$$
Then the quotient group
 $$ \pi_1(R)/N([g_1,h_1], \ [g_2,h_2],\  [g_1,g_2], \ [h_1,g_2],
 \ [h_2^{-1},h_1^{-1}]g_1^{-1},\  [h_2^{-1},g_1^{-1}]g_2^{-1})$$
 surjects to $\pi_1(S)$.

Moreover, the Euler characteristic of $S$, $e(S)$, equals $e(R)+6$ and the signature $\sigma(S)$ equals $\sigma(R)-2$.

 \end{enumerate}

\end{lem}
\begin{proof}  The first assertion is explained in the paragraph that precedes the statement of Lemma \ref{lem1}.   

For the second assertion, the statements about the fundamental group  of $S$  are a straightforward consequence of the Seifert-Van Kampen theorem applied to the decomposition $S = (R-nbd(F))\cup_{\tilde{\phi}}(W-nbd(G))$, using the fact that the meridian of $G$ bounds a disk in $W$ (the punctured exceptional sphere) and that $\pi_1(G)\to \pi_1(W-G)$ is surjective (because its composite with the isomorphism $\pi_1(W-G)\to \pi_1(W)$ is surjective).

The only remaining unverified assertions are the claims about Euler characteristic and signature.
The Euler characteristic  of $S$ is computed using
the formula $e(A\#_H B)=e(A)+e(B)-2e(H)$, which is true for any
  sum of 4-manifolds along surfaces.  Therefore
$e(S)=e(W)+e(R)+4=2+e(R)+4=e(R) +6$.  Novikov additivity can be used
to compute the signature, so
$\sigma(S)=\sigma(W)+\sigma(R)=\sigma(R)-2$.
\end{proof}

In  Lemma \ref{lem1}, suppose further that $R$ is symplectic and $F$ is a symplectic
genus 2 surface in $R$.  Then $S$ admits a symplectic  structure (\cite{Gompf}). Finally, if $R$ is
minimal, and not an $S^2$ bundle over $F$, then $S$ is minimal by Usher's theorem \cite{usher}. This follows
since every embedded $-1$ sphere in $W$ intersects the surface $G$.

\subsection{The manifold $P$}

The second building block $P$ will be the symplectic sum along a torus of two manifolds constructed in the same manner as  $V$.    Alternatively, $P$ can be described as the result of Luttinger surgeries on four Lagrangian tori in the product of a genus two surface with a torus.  There are three perspectives for the reader to keep in mind:

\begin{enumerate}
\item To apply the calculations of Theorem \ref{tough},  one should view $P$ as the union along their boundary of two manifolds obtained by Luttinger surgeries on the product of a punctured torus with a torus, and then apply the Seifert-Van Kampen theorem.
\item To conclude that $P$ is symplectic  one should view $P$ as the symplectic sum of two manifolds obtained by Luttinger surgeries on $T^4=T^2\times T^2$.
\item To conclude that $P$ is minimal one should view $P$ as the symplectic sum of two minimal symplectic manifolds.
\end{enumerate}
Since the fundamental group calculation is the most delicate, we take the first perspective, and trust that the reader can follow the claims about symplectic structure and minimality.

 We therefore build $P$ as the union of two manifolds $P_1$ and $P_2$ along their boundary.  Give each  torus which appears in the following construction the standard symplectic form (i.e.   as the quotient  $\RR^2/\ZZ^2$).  A punctured torus should be given the restricted symplectic form, and the product of two (punctured or unpunctured) tori should be given the product symplectic form.

 \bigskip

 For $P_1$, start with a product $\hat{H}_1\times K_1$ of a torus with base point $h_1$ and a punctured  torus with base point $k_1$.
 Label the loops on $\hat{H}_1$ generating $\pi_1(\hat{H}_1)$ by $x_1,y_1$ and the loops in $K_1$ generating $\pi_1(\hat{K}_1)$ by $s_1,t_1$.   Let $\hat{H}$ and $K$ be as in Theorem \ref{tough} and Scholium \ref{scho}.

 Let $\psi_1:\hat{H}_1 \to \hat{H}$ be the diffeomorphism of the torus which rotates the square by angle $\pi/2$.  Thus $\psi_1$ preserves base points, is orientation preserving,  and induces the isomorphism $x_1\mapsto y$ and $y_1\mapsto x^{-1}$ on fundamental groups. Similarly Let $\psi_2:K_1 \to K$ be the diffeomorphism of the punctured torus which rotates the punctured square by angle $\pi/2$.  Thus $\psi_2$ preserves base points, is orientation preserving,  and induces the isomorphism $s_1\mapsto b$ and $t_1\mapsto a^{-1}$.

Since rotation by $\pi/2$ induces an area-preserving map on the torus, the  diffeomorphism
$\Psi=\psi_1\times\psi_2:  \hat{H}_1\times K_1\to \hat{H}\times K$  is a symplectomorphism which takes the loops $x_1,y_1,s_1,t_1$ to $y,x^{-1},b,a^{-1}$ respectively. We do $+1$ Luttinger surgery on $\Psi^{-1}(T_1)$ (the $y_1^{-1}\times t_1^{-1}$ torus) along $y_1^{-1}$ and
 $+1$ Luttinger surgery on $\Psi^{-1}(T_2)$ (the $x_1\times t_1^{-1}$ torus) along $t_1^{-1}$. Then Theorem \ref{tough} and
 Scholium \ref{scho} imply that the fundamental group of the resulting manifold $P_1$ is generated by
 $$x_1,y_1,s_1,t_1$$ and the Luttinger relations
 $$  y_1=[s_1^{-1}, x_1^{-1}], s_1t_1s_1^{-1}=[y_1,s_1]$$
 as well as the universal relations
 $$[y_1^{-1},t_1^{-1}]=1,\  [x_1,t_1^{-1}]=1,\  [x_1, s_1t_1^{-1}s_1^{-1}]=1,\  [x_1,y_1]=1$$
 hold.   We rewrite the second Luttinger relation as
 $$t_1=[s_1^{-1},y_1].$$

\bigskip

 For $P_2$, start with a product $\hat{H}_2\times K_2$ of a torus and a punctured  torus.
 Label the loops on $\hat{H}_2$ generating $\pi_1(\hat{H}_2)$ by $x_2,y_2$ and the loops in $K_2$ generating $\pi_1(\hat{K}_2)$ by $s_2,t_2$.

As above, choose a symplectomorphism  $\Psi_2:\hat{H}_2\times K_2\to\hat{H}\times K$ which takes the generators $x_2,y_2, s_2,t_2$ to $y,x^{-1}, b, a^{-1}$ respectively.

  We do $+1$ Luttinger surgery on $\Psi_2^{-1}(T_1)$ (the $y_2^{-1}\times t_2^{-1}$ torus) along $t_2^{-1}$ and
 $-1$ Luttinger surgery on $\Psi_2^{-1}(T_2)$ (the $x_2\times t_2^{-1}$ torus) along $x_2$. Then Theorem \ref{tough} and
 Scholium \ref{scho} imply that the fundamental group of the resulting manifold $P_2$ is generated by
 $$x_2,y_2,s_2,t_2$$ and the Luttinger relations
 $$t_2=[s_2^{-1},x_2^{-1}], x_2=[y_2, s_2]$$
  as well as the universal relations
 $$[y_2^{-1},t_2^{-1}]=1, [x_2,t_2^{-1}]=1, [x_2, s_2t_2^{-1}s_2^{-1}]=1,[x_2,y_2]=1$$
 hold.

\bigskip

Denote by $M_1$ and $M_2$ the symplectic manifolds obtained by the same construction as $P_1$ and $P_2$ but starting with closed tori, i.e. $\hat{H}_i\times\hat{K}_i\cong T^4$.  Denote by $z_1$ and $z_2$ the centers of the disks removed from $\hat{K}_i$ to obtain $K_i$.  As a smooth manifold, the symplectic sum, $P$, of $M_1$ and $M_2$ along the symplectic tori with trivial normal bundles
$\hat{H}_1 \times \{z_1\}$ and $\hat{H}_2 \times \{z_2\}$ (\cite{Gompf}), is the union of $P_1$ and $P_2$ along their boundary 3-tori.  We use the diffeomorphism of the tori along which the symplectic sum is performed so that $x_1$ is identified with $x_2$ and $y_1$ is identified with $y_2$.

 More precisely,  there exists an arc 
  $\beta$  in $K_1$ which starts at a point   $k_1'\in\partial K_1$  and ends at $k_1$ and which misses 
  $\psi_2^{-1}(A_i), i=1,2$, since  cutting the surface $K$ along $a\cup b\cup A_1\cup A_2$  does not disconnect $k$ from $\partial K$.   This arc $\beta$ should be (and can be) chosen so that the loop traced out by the boundary is homotopic rel endpoint to $\beta*[s_1,t_1]*\beta^{-1}$ in $K_1$.   The 
  arc $\tilde{\beta}=\{h_1\}\times \beta\subset  \hat{H}_1 \times K_1$ misses $T_1\cup T_2$, since $\beta$ misses $A_1\cup A_2$, and hence can be viewd as a path in $P_1$. 
  
  Conjugating by $\tilde{\beta}$    induces an isomorphism $\pi_1(P_1, (h_1,k_1))\cong 
   \pi_1(P_1, (h_1,k_1'))$ so that the  loops $x_1,y_1,s_1,t_1$ are sent to  loops we temporarily call $x_1',y_1', s_1' , t_1'$.  Obviously, all the relations listed above involving the  $x_1,y_1,s_1,t_1$  also hold for the $x_1',y_1', s_1' , t_1'$.    
   
The loop $x_1'=\tilde{\beta}*x_1*\tilde{\beta}^{-1}$ is homotopic rel endpoint into the boundary of $P_1$. In fact, the one parameter family of loops $x_1'(s)   = \tilde{\beta}_s *(x_1\times \beta(s))*\tilde{\beta}_s^{-1}$  (where $\tilde{\beta}_s(t)=\tilde{\beta}(st)$) gives a homotopy of  $x_1'$  to the loop $x_1\times\{k_1'\}$ in the boundary 3-torus $ \hat{H}_1\times \partial K_1$ of $P_1$. (Note that this uses the fact that $\beta$ misses $\psi_2^{-1}(A_1\cup A_2)$.) A similar comment applies to $y_1'$.  The loop $[s_1', t_1']$ lies entirely on $\{h_1\}\times K_1\subset \hat{H}_1\times K_1-nbd(T_1\cup T_2)\subset P_1$.  Hence the loop $\{h_1\}\times \partial K_1\subset\partial P_1$ maps to $[s_1',t_1']$ via the inclusion $\pi_1(\partial P_1,(h_1,k'))\to \pi_1(P_1,(h_1,k'))$ . 

Thus we  abuse notation slightly and rename $x_1= x_1'(1)$, $y_1=y_1'(1), s_1=s_1' $, and $t_1=t_1'$.
These   loops are based at the base point $(h_1, k_1')$ on the boundary of $P_1$, generate $\pi_1(P_1)$, and all the relations listed above hold. Moreover the three loops $x_1,y_1, $ and $c=\{h_1\}\times \partial K_1$   all lie on the  boundary $\partial P_1$, generate $\pi_1(\partial P_1,(h_1,k_1'))$, and   $c$ is sent to $[s_1,t_1]$ in $\pi_1(P_1)$.

  A similar comment applies to $P_2$, so we end up with the same presentation, but with base point $(h_2,k_2')\in \partial P_2$, and the loops $x_2,y_2, $ and $[s_1,t_1]$ in $\partial P_2$ generating  $\pi_1(\partial P_2,(h_2,k_2'))$

\medskip

 We glue $P_1$ to $P_2$ using a  base point preserving diffeomorphism which takes $x_1 $ to $x_2$,   $y_1$ to $y_2$, and $[s_1,t_1]$ to $[s_2,t_2]^{-1}$.  (This last gluing actually follows from the first two and the fact that we are forming the symplectic sum of  $M_1$ and $M_2$ to build $P$.) Note that we can arrange this to be  the relative  symplectic sum  (\cite{Gompf}) of $(M_1, \{h_1\}\times \hat{K}_1)$ and 
 $(M_2, \{h_2\}\times \hat{K}_2)$ so  that the surfaces $\{h_1\}\times K_1$ and $\{h_2\}\times K_2$  line up along their boundary, yielding a closed symplectic genus 2 surface $F$ in $P$. The loops $s_1,t_1,s_2,t_2$ lie   on $F$ and these form the standard set of generators of the fundamental group  of $F$.   This fact allows us to apply Lemma~\ref{lem1} in the proof of 
Theorem \ref{ex1} below.    The Seifert-Van Kampen theorem implies that $\pi_1(P)$ is generated by
$x_1,y_1,s_1,t_1,x_2,y_2,s_2,t_2$.

\medskip

The definition of $P$,  the  calculations for $P_1$ and $P_2$  given above,  and the Seifert-Van Kampen theorem  imply that the relations
 \begin{equation}\label{rels1}    y_1=[s_1^{-1}, x_1^{-1}], t_1=[s_1^{-1},y_1] , \ [y_1^{-1},t_1^{-1}]=1,\  [x_1,t_1^{-1}]=1,\  [x_1, s_1t_1^{-1}s_1^{-1}]=1,\  [x_1,y_1]=1   \end{equation}

 \begin{equation}\label{rels2}t_2=[s_2^{-1},x_2^{-1}], x_2=[y_2, s_2],
 [y_2^{-1},t_2^{-1}]=1, [x_2,t_2^{-1}]=1, [x_2, s_2t_2^{-1}s_2^{-1}]=1,  [x_2,y_2]=1\end{equation}
and
 \begin{equation}\label{rels3}x_1=x_2,y_1=y_2\end{equation}
hold in $\pi_1(P)$.  The additional relation $[s_1,t_1][s_2,t_2]=1$ also follows from the Seifert-Van Kampen theorem, but we will not need it below.

\medskip

The closed symplectic manifolds $M_1$ and $M_2$ have trivial second homotopy group, and hence are minimal, by Lemma \ref{mini}.   Thus by Usher's theorem \cite{usher}  their symplectic sum $P$ is also minimal.

\bigskip

 The Euler characteristic
is $e(P) = e(M_1) + e(M_2) +0 = 0$ and the signature $\sigma(P) =
\sigma(M_1) + \sigma(M_2) =0$.

  \section{Assembly: an exotic symplectic $\CP^2\#3\overline{\CP}^2$}

 Let $X$ be the symplectic  sum of
$P$ and $W$ along the    genus 2 surfaces $F\subset P$ and $G\subset W$,
$$X=(P-nbd(G))\cup_{\tilde{\phi}}(W-nbd(G))$$
 using a diffeomorphism $\phi:F\to G$ that identifies generators in $\pi_1(F)$ with their namesakes in $\pi_1(G)$.

By  Lemma~\ref{lem1} and the text that immediately follows it, $X$
is a symplectic 4-manifold with $e(X) = 6$ and $\sigma(X)=-2$.
Furthermore, $X$ is minimal by Usher's theorem \cite{usher} since $P$ is, and since    $W-G$ contains no $-1$ spheres.

  Once we show  $X$ is simply
connected, then   Freedman's theorem \cite{Freedman} implies that $X$ is
homeomorphic to $\CP^2\#3\bCP^2$. It  cannot  be diffeomorphic  $\CP^2\#3\bCP^2$
however, since $X$ is minimal, and by  results of Taubes \cite{taubes1,taubes2}, a minimal symplectic 4-manifold  cannot    contain a smoothly
embedded $-1$ sphere, but $\CP^2\#3\bCP^2$ contains smoothly embedded $-1$ spheres, namely,  the exceptional spheres.

\begin{thm}\label{ex1}  The minimal symplectic manifold $X$ is simply connected, hence
homeomorphic but not diffeomorphic to $\CP^2\#3\bCP^2$.

\end{thm}
\begin{proof}

Since the loops $s_1,t_1,s_2,t_2$ lie on $F$, Lemma~\ref{lem1} implies that the fundamental group of $X$ is a quotient of
$\pi_1(P)/N$ where $N$ is the normal subgroup generated by
\begin{equation}\label{rels4} [s_1,t_1],\ [s_2,t_2], \ [s_1,s_2],\ [t_1,s_2],\ [t_2^{-1},t_1^{-1}]s_1^{-1},\   [t_2^{-1},s_1^{-1}]s_2^{-1} .\end{equation}

Denote by {\em relations 1-20} the  14 relations listed for the fundamental group of $P$ in Equations (\ref{rels1}), (\ref{rels2}), and (\ref{rels3}) and the  six additional relations  of Equation (\ref{rels4}).  Recall that $[r,s]^{-1}=[s,r]$.

 To start, observe that relations 1 and 19 imply
 $$y_1=[s_1^{-1},x_1^{-1}]=[[t_1^{-1},t_2^{-1}], x_1^{-1}].$$
Relation 4 implies that $x_1$ commutes with $t_1$ and relations 10 and 13 imply that $x_1$ commutes with $t_2$. This implies that $y_1=1$.

The rest of the generators are rapidly killed. Relation 14 implies  $y_2=1$. Relation 2 implies $t_1=1$.
 Relation 19 implies that $s_1=1$. Relation 20 implies that $s_2=1$. Relation 7 now shows that $t_2=1$ and Relations 8 and 13  imply that $x_1=x_2=1$.

  Thus Lemma \ref{lem1} says that  $\pi_1(P)$ is a quotient of the trivial group, hence is trivial.  As explained above this implies that $X$ is homeomorphic to, but not diffeomorphic to $\CP^2\#3\bCP^2$.
\end{proof}





\bigskip

 
\begin{thebibliography}{99999}

\bibitem{ADK} D. Auroux, S.K. Donaldson, and L. Katzarkov, {\em Luttinger surgery along Lagrangian tori and non-isotopy for singular symplectic plane curves}. Math. Ann. 326 (2003), no. 1, 185--203.

\bibitem{A} A. Akhmedov, {\em Exotic smooth structures on $3\CP^2\#7\bCP^2$}. preprint (2006)
math.GT/0612130

\bibitem{AP} A. Akhmedov, B. Doug Park, {\em Exotic Smooth Structures on small 4-Manifolds.} preprint (2007) math.GT/0701664.

\bibitem{B1} S. Baldridge, {\em New symplectic 4-manifolds with $b_+=1$.} Mathematische Annalen 333 (2005) 633-643.


\bibitem{BK} S. Baldridge and P. Kirk, {\em On symplectic 4-manifolds with prescribed fundamental group.} To appear in
Commentarii Mathematici Helvetici.

\bibitem{BK2} S. Baldridge and P. Kirk, {\em Luttinger surgery and
interesting symplectic 4-manifolds with small Euler characteristic},
preprint. math.GT/0701400.


\bibitem{BK3} S. Baldridge and P. Kirk, {\em Symplectic 4-manifolds with arbitrary fundamental group near the Bogomolov-Miyaoka-Yau line.}
J. Symplectic Geom. 4 (2006), no. 1, 63--70.



\bibitem{D} S. Donaldson, {\em
Irrationality and the $h$-cobordism conjecture. }
J. Differential Geom. 26 (1987), no. 1, 141--168.

 \bibitem{DS} D. McDuff and D. Salamon,
 `Introduction to symplectic topology.
Second edition.'  Oxford Mathematical Monographs. The Clarendon Press, Oxford University Press, New York, 1998. x+486

\bibitem{FS3} R. Fintushel and R. Stern, {\em Double node neighborhoods and families of simply connected 4-manifolds with $ b^+=1$.}  J. Amer. Math. Soc. 19 (2006), no. 1, 171--180.

\bibitem{FS4} R. Fintushel and R. Stern, {\em
Families of simply connected 4-manifolds with the same Seiberg-Witten invariants.}
Topology 43 (2004), no. 6, 1449--1467.


\bibitem{Freedman} M. Freedman, {\em The topology of four-dimensional  manifolds.}
J. Differential Geom. 17 (1982), no. 3, 357--453.

\bibitem{Gompf} R. Gompf, {\em A new construction of symplectic manifolds}.  Ann. of Math. (2) 142 (1995), no. 3, 527--595.


\bibitem{Kot} D. Kotschick, {\em  On manifolds homeomorphic to $\CP^2\#8\bCP^2$.}  Invent. Math. 95 (1989), no. 3, 591--600.

\bibitem{Lut} K. M. Luttinger, {\em Lagrangian Tori in $\RR^4$}. J. Diff. Geom.  {52} (1999), 203--222.


\bibitem{OzS} P. Ozsv\'{a}th and Z. Szab\'{o}, {\em On Park's exotic smooth four-manifolds,} `Geometry and topology of manifolds', 253--260, Fields Inst. Commun., { 47}, Amer. Math. Soc., Providence, RI, 2005.


\bibitem{park1} Jongil Park, {\em
Simply connected symplectic 4-manifolds with $b\sp +\sb 2=1$ and $c\sp 2\sb 1=2$. }
Invent. Math. 159 (2005), no. 3, 657--667.


\bibitem{PSS} Jongil Park, A.  Stipsicz, and Z. Szabo, {\em
Exotic smooth structures on $\Bbb{CP}\sp 2\#5\overline{\Bbb{CP}\sp 2}$.}
Math. Res. Lett. 12 (2005), no. 5-6, 701--712.


\bibitem{SS} A. Stipsicz and Z. Szab\'{o}, {\em An exotic smooth structure on $\CP^2\#6\overline{\CP}^2$,} Geom. Topol. {9} (2005), 813--832.


\bibitem{taubes1} C. Taubes, {\em Counting pseudo-holomorphic submanifolds in dimension  4.} J. Differential Geom. 44 (1996), no. 4, 818--893.

\bibitem{taubes2} C. Taubes, {\em Seiberg-Witten and Gromov invariants.} Geometry and physics (Aarhus, 1995), 591--601, Lecture Notes in Pure and Appl. Math., 184, Dekker, New York, 1997.


\bibitem{usher} M. Usher, {\em Minimality and symplectic sums}. To appear in Internat. Math. Res. Not.  preprint (2006) math.SG/0606543

\end{thebibliography}
\end{document}